\documentclass{article}

\usepackage{amsmath}
\usepackage[francais]{babel}
\usepackage{amssymb}

\def\ni{\noindent}
\def\Q{\hbox{\bf Q}}
\def\F{\hbox{\bf F}}
\def\P{\hbox{\bf P}}
\def\Z{\hbox{\bf Z}}
\def\C{\hbox{\bf C}}
\def\G{\hbox{\bf G}}
\def\Aut{\mathop{\rm Aut}}
\def\Weyl{\mathop{\rm Weyl}}
\def\deg{\mathop{\rm deg}}
\def\det{\mathop{\rm det}}
\def\GL{\mathop{\bf GL}}
\def\SL{\mathop{\bf SL}}
\def\PGL{\mathop{\bf PGL}}
\def\PSL{\mathop{\bf PSL}}
\def\End{\mathop{\rm End}}
\def\Gal{\mathop{\rm Gal}}
\def\charb{\mathop{\rm char}}
\def\Cr{\mathop{\rm Cr}}
\def\Pic{\mathop{\rm Pic}}
\def\Aut{\mathop{\rm Aut}}
\def\Hom{\mathop{\rm Hom}}

\def\l{\ell}

\def\iif{{\rm if}}
\def\oor{{\rm or}}

\begin{document}

\begin{center}
    {\sc A Minkowski-style bound for the orders of the finite subgroups
                      of the Cremona group of rank $2$ over an arbitrary field} 
\end{center}

\medskip

\begin{center}

                        {\bf Jean-Pierre Serre}
\end{center}

\medskip

\begin{center}

\end{center}

\bigskip

   Let  $k$  be a field. Let  $\Cr(k)$  be the Cremona group of rank $2$
over  $k$, i.e. the group of $k$-automorphisms of  $k(X,Y)$, where  $X$  and  $Y$
are two indeterminates.

   We shall be interested in the finite subgroups of  $\Cr(k)$  of
order prime to the characteristic of  $k$. The 
case  $k = \C$  has a long history, going back to the $19$-th century
(see the references in  [Bl 06] and [DI 07]), and culminating in
an essentially complete (but rather complicated) classification, see  [DI 07].
For an arbitrary field, it seems reasonable to simplify the problem
\`a la Minkowski, as was done in  [Se 07]  for semisimple groups; this 
means giving a sharp multiplicative bound for the orders of the 
finite subgroups we are considering.

   In \S  $6.9$ of [Se 07], one finds a few questions in that direction,
for instance the following:

   If  $k = \Q$, is it true that  $\Cr(k)$  does not contain any element of prime
order $\geqslant 11$ ?

  More generally, what are the prime numbers  $\l$, distinct from 
$\charb(k)$, such that  $\Cr(k)$  contains an element of order  $\l$ ? \medskip

  This question has now been solved by Dolgachev and Iskovskikh ([DI 08]), the answer being that  there is equivalence between:

\medskip

     $\Cr(k)$  contains an element of order  $\l$

\noindent
and

    $[k(z_{\l}):k] = 1,2,3,4$ or $6$, where  $z_{\l}$  is a primitive $\l$-th root of unity.

\medskip

  As we shall see, a similar method can handle arbitrary  $\l$-groups and
one obtains an explicit value for the Minkowski bound of  $\Cr(k)$, in terms
of the size of the Galois group of the cyclotomic extensions of  $k$
(cf. Theorem 2.1 below). For instance:

\medskip
\ni
{\bf Theorem} - {\it Assume  $k$  is finitely generated over its prime subfield. Then
the finite subgroups of  $\Cr(k)$  of order prime to $\charb(k)$ have bounded
order. Let  $M(k)$  be the least common multiple of their orders}.

\smallskip

\ni  {\it  a) If  $k = \Q$, we have  $M(k) = 120960 = 2^7.3^3.5.7.$

\smallskip

\ni b) If  $k$  is finite with  $q$  elements, we have} :

$M(k) = \left\{
\begin{array}{ll}

3.(q^4 - 1)(q^6 - 1) & {\it if } \ q \equiv 4 \ {\it or} \;7 \pmod 9\\

(q^4 - 1)(q^6 - 1)&{\it otherwise}.
\end{array}
\right.
$

\medskip

  For more general statements, see  \S  $2$.  These statements involve the  cyclotomic 
invariants of  $k$  introduced in [Se 07, \S  $6$]; their 
definition is recalled in \S 1. The proofs are given in \S  $3$ (existence 
of large subgroups) and in \S  $4$ (upper bounds). For the upper bounds, 
we use a method introduced by Manin ([Ma 66])
and perfected by Iskovskikh ([Is 79], [Is 96]) and Dolgachev-Iskovskikh
([DI 08]); it allows us to realize any finite subgroup of  $\Cr(k)$
as a subgroup of  $\Aut(S)$, where  $S$  is either a Del Pezzo surface or 
a conic bundle over a conic.
  A few conjugacy  results are given in  \S  $5$. The last $\S$ contains a series of questions on the Cremona groups of rank $> 2.$

\bigskip

\section{The cyclotomic invariants  $t$  and  $m$}
In what follows, $k$  is a field,  $k_s$  is a separable closure of  $k$
and  $\overline{k}$ is the algebraic closure of  $k_s$.

   Let  $\l$ be
a prime number distinct from  $\charb(k)$; the $\l$-adic valuation of  $\Q$ 
is denoted by  $v_{\l}$.   If  $A$  is a finite set, with cardinal $|A|$,  we write
   $v_{\l}(A)$  instead of  $v_{\l}(|A|)$.

   There are two invariants   $t = t(k,l)$  and  $m = m(k,l)$  
which are associated with the pair  $(k,l)$, cf.
[Se 07, \S 4]. Recall their definitions:

\subsection{Definition of  $t$}

      Let  $z \in k_s$  be
  a primitive $\ell$-th root of unity  if  $l > 2$  and  a primitive 
$4$-th root of unity
if  $\ell = 2$. We put  $$t = [k(z):k].$$  If  $\l > 2$,  $t$  divides  $\l-1$. If
$\l = 2$ or $3$, then  $t = 1$ or $2$.

\subsection{Definition of  $m$}

For  $\l > 2$,  $m$  is the upper bound (possibly infinite)  of the  $n$'s
  such that $k(z)$  contains the $\l^n$-th roots of unity. We have  $m \geqslant 1$.

   For  $\l = 2$, $m$  is the upper bound (possibly infinite) of the  $n$'s
such that  $k$  contains  $z(n) + z(n)^{-1}$ , where  $z(n)$  is a primitive
$2^n$-root of unity. We have  $m \geqslant 2$. [The definition of  $m$  given 
in  [Se 07, $\S4.2$]  looks different, but it is equivalent to the one 
here.]

\smallskip
\ni
{\it Remark}. When $\l > 2$, knowing  $t$ and $m$  amounts to knowing the image of the
$\l$-th cyclotomic character $\Gal(k_s/k)\rightarrow \Z_l^*$,  cf. [Se 07, \S 4].

\subsection{Example : $k=\Q$}
Here, $t$  takes its largest possible value, namely  $t = \l-1$  for
     $\l > 2$ and  $t = 2$  for  $\l = 2$. And  $m$  takes its smallest 
possible
     value, namely  $m = 1$  for  $\l > 2$  and  $m = 2$  for  $\l = 2$.

\subsection{Example : $k$  finite with  $q$  elements}
      If $\l > 2$, one has:

\begin{tabular}{l}

        $t$  =  order of  $q$  in the multiplicative group  $\F^*_{\l}$\\
        
        \smallskip
        
        $m$  =  $v_{\l}(q^t - 1) = v_{\l}(q^{l-1} - 1)$.
\end{tabular}

  \ni If $\l = 2$, one has:
   
  \begin{tabular}{l}
  
    $t$ = order of $q$ in $(\Z/4 \Z)^*$\\
    $m = v_2(q^2-1) - 1$.
    \end{tabular}

\section{Statement of the main theorem}

  Let  $K = k(X,Y)$, where  $X,Y$  are indeterminates, and let  $\Cr(k)$  be the
Cremona group of rank $2$ over  $k$, i.e. the group  $\Aut_k K$ .
   Let  $\l$  be a prime number, distinct from $\charb(k)$,
and let  $t$  and  $m$
be the cyclotomic invariants defined above.

\subsection{Notation}
Define a number  $M(k,\l) \in \{0,1,2,...,\infty\}$ as follows:

\medskip

      For  $\l = 2,  \;\;M(k,\l) = 2m + 3$.

\medskip

$
{\rm For}\;\;  \l = 3,  \;\; M(k,\l) = 
\left\{
\begin{array}{ll}

4&\iif  \;\;t = m = 1\\

                     2m + 1&{\rm otherwise}.

\end{array}
\right.
$

\medskip

$
{\rm For}\;\;  \l > 3,  \;\;M(k,\l)=
\left\{
\begin{array}{lll}

  2m\;\;&  \iif  \;\;t = 1 \; \; \oor \; \, 2\\

                 m\;\;&   \iif  \;\;t = 3, 4 \; \;\oor \; \, 6\\

                  0\;\;&   \iif  \;\;t = 5\; \;  \oor \; \, t > 6.
\end{array}
\right.
$

\subsection{The main theorem} 

\medskip

\ni
{\bf Theorem  $2.1$}.(i)  {\it Let  $A$  be a finite subgroup of  $\Cr(k)$. Then
$v_{\l}(A)  \leqslant M(k,\l)$.

  {\rm (ii)} Conversely, if  $n$  is any integer $\geqslant 0$  which is  $\leqslant  M(k,\l)$, then
$\Cr(k)$  contains a subgroup of order  $\l^n$.

 {\rm (In other words, $M(k,\l)$  is the upper bound of the  $v_{\l}(A).$)}} 

\medskip

The special case where  $A$  is cyclic of order  $\l$  gives:

\medskip

\ni
{\bf Corollary $2.2$} ([DI 08]). {\it The following properties are equivalent} :

   {\rm a)} {\it $\Cr(k)$  contains an element of order  $\l$

   {\rm b)}  $\varphi(t)\leqslant 2$, i.e. $t = 1, 2, 3, 4$ or $6$.
}

\medskip

   Indeed, b)  is equivalent to  $M(k,\l) > 0$.

\subsection{Small fields}
Let us say that  $k$  is {\it small} if it has the following properties:

\smallskip

\begin{tabular}{ll}

     (2.3.1)&    $m(k,\l) < \infty$ for every  $\l \neq \charb(k)$\\
     (2.3.2) &   $t(k,\l)  \rightarrow  \infty $ when  $\l \rightarrow \infty$.
\end{tabular}

\medskip

\ni
{\bf Proposition $2.3$}. {\it A field which is finitely generated over  $\Q$  or  $\F_p$
is small.}
\medskip

{\it Proof}. The formulae given in $\S1.3$ and $\S1.4$ show that both  $\F_p$  and  $\Q$  are
small.  If  $k'/k$  is a finite extension, one has
$$[k':k].t(k',\l) \geqslant  t(k,\l) \;\;   {\rm and}\;\;   m(k',\l) \leqslant m(k,\l) + 
\log_{\l}([k':k]),$$
which shows that  $k$  small $\Rightarrow k'$  small. If  $k'$  is a regular extension
of  $k$ , then  
$$t(k',\l) = t(k,\l) \;\;\;{\rm and} \;\;\; m'(k',\l) = m(k,\l),$$ 
which also shows that
$k$ small $\Rightarrow k'$ small. The proposition follows.

\medskip

 Assume now that  $k$  is small. We may then define an integer  $M(k)$
by the following
formula

\medskip

   $(2.3.3)\;\;\;\;\;\;\;\;\;\;\;\;\;\;\;\;\;
M(k) =  \displaystyle{\prod_{\l} \l^{M(k,\l)}}$,

\medskip

\ni
where  $\l$  runs through the prime numbers distinct from  $\charb(k$).
The
formula makes sense since  $M(k,\l)$ is finite for every  $\l$  and is  $0$
for every   $\l$  but a finite number. With this notation, Theorem $2.1$ can be
reformulated as:
\medskip

\ni
{\bf Theorem $2.4$}. {\it If  $k$  is small, then the finite subgroups of  $\Cr(k)$
of order prime to $\charb(k)$ have bounded order, and the l.c.m. of their
orders is the integer  $M(k)$  defined above.}

\medskip
Note that this applies in particular when  $k$  is finitely generated
over its prime subfield.

\subsection{
Example: the case $k = \Q$}

By combining $1.3$ and $2.1$, one gets

$$
M(\Q,\l) = 
\left\{
\begin{array}{lll}
7&  {\rm for}&  \l = 2,\\
3&  {\rm for} &  \l = 3,\\
1&  {\rm for} &  \l = 5,7\\
0&  {\rm for} &  \l > 7.
\end{array}
\right.
$$

This can be summed up by:
\medskip

\ni
{\bf Theorem  $2.5$}. \  {\it    $M(\Q) = 2^7.3^3.5.7$}.
\medskip

\subsection{Example : the case of a finite field}

\ni
{\bf Theorem $2.6$}. {\it If  $k$  is a finite field with  $q$  elements, we have
$$M(k) =
\left\{
\begin{array}{ll}
     3.(q^4 - 1)(q^6 - 1)&   if  \;\;q \equiv 4 \;\ or \; 7 \pmod 9\\
 \smallskip
         
     (q^4 - 1)(q^6 - 1)  & otherwise.
\end{array}
\right.
$$
}

\medskip

{\it Proof}. Denote by  $M'(k,\l)$  the $\l$-adic valuation of the right side of
the formulae above.

   If  $\l$  is not equal to  $3$,  $M'(k,\l)$  is equal to
      $$v_{\l}(q^4 - 1) + v_{\l}(q^6 - 1)$$
and we have to check that  $M'(k,\l)$  is equal to  $M(k,\l).$

   Consider first the case  $\l = 2$. It follows from the definition 
of  $m$  that

  $v_2(q^2 - 1) = m + 1$, and hence  $v_2(q^4 - 1) = m + 2$  and  $v_2(q^6 - 1) =
m + 1$. This gives  $M'(k,\l) = 2m + 3 = M(k,\l).$

   If  $\l > 3$, the invariant  $t$  is the smallest integer $> 0$  such that
$q^t = 1 \pmod \l$. If  $t = 5$  or $t > 6$, this shows that  $M'(k,\l) = 0$.

If  $t = 3$  or  $6$, $q^4 - 1$  is not divisible by  $\l$  and  $q^6 - 1$  is
divisible by  $\l$; moreover, one has  $v_{\l}(q^6 - 1) = m$. This gives
$M'(k,\l) = m = M(k,\l)$. Similarly, when  $t = 4$, the only factor 
divisible by $\l$  is  $q^4 - 1$  and its $\l$-adic valuation is  $m$. 
When  $t = 1 $ or $2$, both factors
are divisible by  $\l$  and their $\l$-adic valuation is  $m$.

   The argument for  $\l = 3$  is similar:  we have 
$$v_3(q^4 - 1)=m\;\;\;{\rm and} \;\;\;
v_3(q^6 - 1)=m + 1.$$
 The congruence  $q \equiv 4$ 
or $7 \pmod 9$  means that  $t = m = 1$.

\medskip

For instance:

\smallskip

\ni
$M(\F_2) = 3^3.5.7; \ \ \   M(\F_3) = 2^7.5.7.13;\ \ \   M(\F_4) = 3^4.5^2.7.13.17;$

\smallskip

$ \ \ M(\F_5) = 2^7.3^3.7.13.31; \  M(\F_7) = 2^9.3^4.5^2.19.43.$

\subsection{Example: the $p$-adic field  $\Q_p$}

For  $\l \neq p$, the  $t,m$  invariants of  $\Q_p$  are the same as those
    of  $\F_\l$ , and for  $\l = p$  they are the same as those of  $\Q$. 

This
shows that  $\Q_p$  is `` small '', and a simple computation gives

         $$M(\Q_p) = c(p).(p^4 - 1)(p^6 - 1),$$
with

$c(2) = 2^7; \ c(3) = 3^3; \ c(5) = 5; \ c(7) = 3.7;$

      $c(p) = 3$ \  if \  $p > 7$ \  and   \ $p \equiv 4$ or $7 \pmod 9;$

      $c(p) = 1$\  otherwise.

\medskip

For instance: 

\smallskip

\ni
$M(\Q_2) = 2^7.3^3.5.7; \ \ M(\Q_3) = 2^7.3^3.5.7.13;\ \
M(\Q_5) = 2^7.3^3.5.7.13.31; $ 

\smallskip

\ni
\ \ \ \ \ \  $M(\Q_7) = 2^9.3^4.5^2.7.19.43;\ \
M(\Q_{11}) = 2^7.3^3.5^2.7.19.37.61.$

\subsection{Remarks}

   1. The statement of Theorem 2.6 is reminiscent of the formula which gives
the order of  $G(k)$, where  $G$  is a split semisimple group and $|k| = q$.  In 
such a formula,
the factors have the shape  $(q^d - 1)$, where  $d$  is an invariant degree
of the Weyl group, and  the number of factors is equal to the rank of  $G$.
Here also the number of factors is equal to the rank of $\Cr$,  which is  $2$.
The exponents  $4$  and $6$  are less easy to interpret. In the proofs below, they occur as the maximal orders of the torsion elements
of the `` Weyl group '' of  $\Cr$, which is  $\GL_2(\Z)$. See also $\S6$.

\smallskip
\noindent
  2. Even though Theorem 2.6 is a very special case of Theorem 2.1, it contains
almost as much information as the general case. More precisely, we could
deduce Theorem 2.1.(i) [which is the hard part] from Theorem 2.6 by the Minkowski
method of reduction$\pmod p$ explained in [Se $07$, \S 6.5].

\smallskip
\noindent
  3. In the opposite direction, if we know Theorem 2.1.(i) for fields of
  characteristic $0$ (in the slightly more precise form given in \S 4.1),
  we can get it for fields of characteristic $p > 0$ by lifting over the
  ring of Witt vectors; this is possible: all the cohomological obstructions  vanish
  (for a detailed proof, see [Se 08, \S5]).

\smallskip
\noindent
  4. For large fields, the invariant  $m$  can be $\infty$. If  $t$  is not
$1,2,3,4$ or $6$, Corollary 2.2 tells us that  $\Cr(k)$ is $\l$-torsion-free. But if  $t$
is one of these five numbers, the above theorems tell us nothing. Still,
as in [Se $07$, \S 14, Theorem 12 and Theorem 13]  one can prove the following:

   a) If  $t = 3,4$ or $6$, then  $\Cr(k)$  contains a subgroup isomorphic
to  $\Q_{\l}/\Z_{\l}$  and does not contain  $\Q_{\l}/\Z_{\l} \times \Q_{\l}/\Z_{\l}$.

   b) If  $t = 1$ or $2$, then  $\Cr(k)$  contains a subgroup isomorphic to
$\Q_{\l}/\Z_{\l} \times \Q_{\l}/\Z_{\l}$ and does not contain a product of three copies of  $\Q_{\l}/\Z_{\l}$.

\section{Proof of Theorem 2.1.(ii)}
   We have to construct large $\l$-subgroups of  $\Cr(k)$. It turns out
that we only need two constructions, one for the very special case
$\l = 3, t = 1, m=1$, and one for all the other cases.

\subsection{The special case  $\l = 3, \ t = 1,\ m=1$}

We need to construct a subgroup of  $\Cr(k)$ of order  $3^4$. To do so we
use the Fermat cubic surface  $S$  given by the homogeneous equation

      $$x^3 + y^3 + z^3 + t^3 = 0.$$

It is a smooth surface, since  $p \neq 3$. The fact that  $t = 1$  means that
$k$  contains a primitive cubic root of unity. This implies that 
the $27$ lines
of  $S$  are defined over  $k$, and hence  $S$  is $k$-rational: its function
field is isomorphic to  $K = k(X,Y)$.  Let  $A$  be the group of automorphisms
  of $S$  generated by the two elements

   $$(x,y,z,t) \mapsto (rx,y,z,t)\;\;\;{\rm and}\;\;\; (x,y,z,t)  \mapsto (y,z,x,t)$$
  where  $r$  is a primitive $3$-rd root of unity.

  We have  $|A| = 3^4$  and  $A$  is a subgroup of  $\Aut(S)$, hence a 
subgroup of  $\Cr(k)$.

\subsection{The generic case}

   Here is the general construction:

   One starts with a 2-dimensional torus  $T$  over  $k$  , with an $\l$-group  $C$\textbf{}
acting faithfully on it. Let  $B$  be an $\l$-subgroup of  $T(k)$. Assume that
$B$  is stable under  $C$, and let  $A$  be the semi-direct product  $A = B.C$.
If we make  $B$  act on the variety  $T$  by translations, we get an action
of  $A$, which is faithful. This gives an embedding of  $A$  in $\Aut(k(T))$,
where  $k(T)$  is the function field of  $T$. By a theorem of Voskresinski\u{\i}
(see [Vo $98$, \S 4.9])
$k(T)$  is isomorphic to  $K = k(X,Y)$. We thus get an embedding of  $A$
in  $\Cr(k)$. Note that  $B$ is {\it toral}, i.e. is contained in the $k$-rational points of a maximal torus of  $\Cr$.

   It remains to explain how to choose $T$, $B$ and $C$. We shall define  $T$
by giving the action of  $\Gamma_k = \Gal(k_s/k)$  on its character
group ; this amounts to giving an homomorphism  $\Gamma_k  \rightarrow \GL_2(\Z)$.

\medskip

\subsubsection{The case  $\l = 2$}

 Let  $n$  be an integer  $\leqslant  m$. If  $z(n)$  is a primitive $2^n$-root
of unity, $k$  contains  $z(n) + z(n)^{-1}$. The field extension  $k(z(n))/k$
has degree  $1$  or  $2$, hence defines a character  $\Gamma_k  \rightarrow  \{1,-1\}$.
Let  $T_1$  be the $1$-dimensional torus associated with this character.
If  $k(z(n)) = k$ , $T_1$ is the split torus  $\mathbf{G}_m$  and we have  $T_1(k) = k^*$.
If  $k(z(n))$ is quadratic over  $k$, $T_1(k)$  is the subgroup of  $k(z(n))^*$
made up of the elements of norm $ 1$. In both cases, $T_1(k)$  contains  $z(n)$.
We now take for  $T$  the torus  $T_1 \times T_1$ and for  $B$  the subgroup of
elements of  $T$  of order dividing  $2^n$. We have  $v_2(B) = 2n$. We take
for  $C$  the group of automorphisms generated by  $(x,y)  \mapsto  (x^{-1},y)$
and  $(x,y)  \mapsto  (y,x)$; the group  $C$  is isomorphic to the dihedral
group  $D_4$; its order is  $8$. We then have  $v_2(A) = v_2(B) + v_2(C) = 2n + 3$,
as wanted.

\medskip

(Alternate construction: the group  $\Cr_1(k) = \PGL_2(k)$  contains a dihedral
subgroup  $D$  of order  $2^{n+1}$; by using the natural embedding of
$(\Cr_1(k) \times \Cr_1(k)).2$  in  $\Cr(k)$ we obtain a subgroup of  $\Cr(k)$ 
isomorphic       to $(D \times D).2$, hence of order  $2^{2n+3}$.)

\medskip

\subsubsection{The case  $\l > 2$}

   We start similarly with an integer  $n \leqslant  m$. We may assume that
the invariant  $t$  is equal to  $1,2,3,4$ or $6$; if not we could take  $A = 1$.
Call  $C_t$  the Galois group of  $k(z)/k$ , cf.\S 1. It is a cyclic group
of order  $t$. Choose an embedding of  $C_t$  in  $\GL_2(\Z)$, with the condition
that, if  $t = 2$, then the image of  $C_t$  is  $\{1,-1\}$. The composition map

 \ \ \ \ \ \ \ \ \ \ \    $r :  \Gamma_k  \rightarrow  \Gal(k(z)/k) = C_t  \rightarrow  \GL_2(\Z)$

\ni defines a $2$-dimensional torus  $T$.

   The group  $B$  is the subgroup $T(k)[\l^n]$ of  $T(k)$  made up of elements
of order dividing  $\l^n$ . We take  $C$  equal to  $1$, except when  $\l = 3$
where we choose it of order  $3$  (this is possible since  $t = 1$  or  $2$
for $\l = 3$, and the group of  $k$-automorphisms of  $T$  is isomorphic to
$\GL_2(\Z))$. We thus have:

    $v_\l(A) =  v_\l(B)$  if $\l > 3$    and  $v_\l(A) = 1 + v_\l(B)$  if  $\l = 3$.
    
    \smallskip

   It remains to estimate $v_\l(B)$. Namely:

\medskip

  ($3.2.3$) $v_\l(B)  = 2n$  \ if \ $t = 1$ or $2$

   This is clear if  $t = 1$  because in that case  $T$  is a split torus
of dimension $2$, and  $k$   contains  $z(n)$.

   If $t = 2$, then  $T = T_1 \times T_1$,  where  $T_1$  is associated with
the quadratic character  $\Gamma_k  \rightarrow  \Gal(k(z)/k)$. We may identify
$T_1(k)$  with the elements of norm $1$ of  $k(z)$, and this shows that
$z(n)$ is an element of  $T_1(k)$  of order  $2^n$. We thus get  $v_\l(B) = 2n$.

\medskip

  ($3.2.4)$  $v_\l(|B|) \geqslant  n$ \ if  \ $t = 3, 4$ or $6$

   We use the description of  $T$  given in  [Se $07$, \S 5.3]: let  $L$  be the
field  $k(z)$. It is a cyclic extension of  $k$  of degree  $t$. Let  $s$  be
a generator of  $C_t = \Gal(L/k)$. Let  $T_L = R_{L/k}(\G_m)$  be the torus
`` multiplicative group of  $L$ ''; we have  $\dim T_L = t$, and  $s$  acts on
$T_L$. We have  $s^t - 1 = 0$  in  $\End(T_L)$. Let  $F(X)$  be the cyclotomic
polynomial of index  $t$, i.e.

\medskip

\begin{tabular}{ll}
     $F(X) = X^2 + X + 1$&   if  $t = 3$\\
     $F(X) = X^2 + 1$ &       if  $t = 4$\\
     $F(X) = X^2 - X + 1$&   if  $t = 6$.
\end{tabular}

\medskip

This polynomial divides  $X^t - 1$; let  $G(X)$  be the quotient  $(X^t-1)/F(X)$,
and let  u  be the endomorphism of  $T_1$  defined by  $u = G(s)$. One
checks ({\it loc.cit.}) that the image  $T$  of  $u : T_1  \rightarrow  T_1$  is a 
2-dimensional torus, and  $s$  defines an automorphism   $s_T$  of    $T$  of
order  $t$, satisfying the equation  $F(s_T) = 0$. This shows that  $T$  is
the same as the torus also called  $T$  above. Moreover, it is easy to
check that the element  $z(n)$  of  $T_1(k)$  is sent by  $u$  into an element
of  $T(k)$  of order  $\l^n.$  This shows that  $v_\l(B)  \geqslant n$.

   [When  $t = 3$, we could have defined  $T$  as the kernel of the norm
map  \\ $N : T_1  \rightarrow  \G_m$. There is a similar definition for  $t = 4$, but
the case  $t = 6$  is less easy to describe concretely.]

   This concludes the proof of the ``existence part''  of Theorem 2.1.

\section{Proof of Theorem 2.1.(i)}

\subsection{Generalization}
   In Theorem 2.1.(i), the hypothesis made on the $\l$-group  $A$  is that it is
contained in  $\Cr(k)$. This is equivalent to saying that $A$ is contained in  $\Aut(S)$,  where  $S$
is a $k$-rational surface, cf. 
e.g. [DI $07$, Lemma 6]. We now want to relax this hypothesis: we 
will merely assume that  $S$  is a surface which is `` geometrically 
rational '', i.e. becomes rational over $\overline{k}$;
for instance  $S$
can be any smooth cubic surface in  $\P_3$. In other words, we will 
be interested
in field extensions  $L$  of  $k$  with the property:

(4.1.1)  $\;\;\;\;\;\;\;\overline{k}\otimes L$  is $\overline{k}$-isomorphic to  $\overline{k}(X,Y)$.

   We shall say that a group  $A$  has `` property  $\Cr_k$ '' if it can be
embedded in  $\Aut(L)$, for some  $L$  having property  ($4.1.1$). The bound given
in Theorem 2.1.(i)  is valid for such groups. More precisely:

\medskip

\ni
{\bf Theorem 4.1}.  {\it If a finite  $\l$-group  $A$  has property  $\Cr_k$ , then
   $v_\l(A)  \leqslant  M(k,\l)$,  where  $M(k,\l)$ is as in \S $2.1$. }

\medskip

   This is what we shall prove. Note that we may assume that  $k$  is
perfect since replacing  $k$  by its perfect closure does not change
the invariants  $t,m$ and $M(k,\l)$.

   [As mentioned in \S $2.7$, we could also assume that  $k$  is finite, or, if we preferred to,
that  $\charb(k)=0$. Unfortunately, none of these reductions
is really helpful.]

\subsection{Reduction to special cases}

We start from an $\l$-group  $A$  having property  $\Cr_k$. As explained
above, this means that we can embed
$A$  in  $\Aut(S)$, where  $S$  is a smooth projective $k$-surface, which is
geometrically rational. Now, the basic
tool is the `` minimal model theorem '' (proved in [DI 07, \S $2$]) which
allows us to assume that  $S$  is of one of the following two types:

\medskip

  a) ({\it conic bundle case}) There is a morphism  $f: S  \rightarrow  C$, where 
$C$  is a smooth
genus zero curve, such that the generic fiber of  $f$  is a smooth curve
of genus  $0$. Moreover, $A$  acts on  $C$  and  $f$  is compatible with 
that action.

  b) ({\it Del Pezzo})  $S$  is a Del Pezzo surface, i.e. its anticanonical class
  $- K_S$  is ample.
\medskip

   In case  b), the degree  $\deg(S)$  is defined as  $K_S.K_S$ 
(self-intersection);
one has  $1  \leqslant  \deg(S)  \leqslant 9$.

   We shall look successively at these different cases.
In the 
second case, we shall use without further reference the standard 
properties of the Del Pezzo surfaces; one can find them for 
instance in [De 80], [Do 07], [DI 07], [Ko 96], [Ma 66] and [Ma 86].

\medskip

\ni
{\it Remark}. In some of these references, the ground field is assumed to 
be of characteristic $0$, but there is very little difference in 
characteristic $p > 0$; moreover, as pointed out above, the characteristic $0$ case 
implies the characteristic $p$ case, thanks to the fact that  $|A|$  is 
prime to $\charb(k)$.

\subsection{The conic bundle case}
Let  $f: S \rightarrow  C$  be as in a) above, and let
$A_o$  be the subgroup of  $\Aut(C)$  given by the action of  $A$  on  $C$.
The group  $\Aut(C)$  is a $k$-form of  $\PGL_2$. By using (for instance) 
[Se $07$, Theorem 5] we get :

$$ v_\l(A_o)\leqslant
\left\{
\begin{array}{lll}

   m + 1  &{\rm if}&  \l = 2,\\
   m  &{\rm if}&  \l > 2  \;\;{\rm and}\;\;  t = 1  \;\;{\rm or}\;\;  2,\\
   0  &{\rm if}&  t > 2.
\end{array}
\right.
$$

   Let  $B$  be the kernel of  $A \rightarrow A_o$. The group  $B$  is a subgroup
of the group of automorphisms of the generic fiber of  $f$. This fiber is
a genus $0$ curve over the function field  $k_C$  of  $C$. Since  $k_C$  is a
regular extension of  $k$, the  $t$  and  $m$  invariants of  $k_C$  are the same
as those of  $k$. We then get for  $v_\l(B)$  the same bounds  as for  $v_\l(A_o)$,
and by adding up this gives:

$$v_\l(A)\leqslant 
\left\{
\begin{array}{lll}

       2m + 2&   {\rm if}&  \l = 2\\
      2m&  {\rm if}&  \l > 2  \;\;{\rm and}\;\;  t = 1 \;\;{\rm or}\;\; 2\\
      0&   {\rm if}&\; t > 2.
\end{array}
\right.
$$
    In each case, this gives a bound which is at most equal to the number
$M(k,\l)$  defined in \S 2.1.

\subsection{The Del Pezzo  case : degree $9$}

   Here  $S$  is $\overline{k}$-isomorphic to the projective plane  $\P_2$ ; in other
words, $S$ is a Severi-Brauer variety of dimension $2$. The group  $\Aut 
S$ is an inner $k$-form of  $\PGL_3$. By using  [Se $07$, \S 6.2]  one finds:
$$v_\l(A)\leqslant 
\left\{
\begin{array}{lll}

      2m + 1&  \iif & \l = 2\\
         2m + 1&  \iif&  \l = 3, \;  t = 1\\
         \leq  m + 1   &\iif  &\l = 3,\;  t = 2\\
         \leq  2m  &    \iif&  \l > 3, \; t = 1\\
         \leq   m    &  \iif  &\l > 3, \; t = 2 \;\;{\rm or\;\;} 3\\
          =  0     &  \iif& t > 3.
\end{array}
\right.
$$

   Here again, these bounds are $\leqslant M(k,\l)$.

\subsection{The Del Pezzo case : degree $8$}

    This case splits into two subcases:

   a)  $S$ is the blow up of  $\P_2$  at one rational point. In that case  $A$
acts faithfully on  $\P_2$  and we apply 4.4.

   b)  $S$  is a smooth quadric of  $\P_3$. The connected component  $\Aut^o(S)$
of  $\Aut(S)$  has index $2$. It is a $k$-form of  $\PGL_2 \times \PGL_2$. If we denote
by  $A_o$  the intersection of  $A$  with  $\Aut^o(S)$, we obtain, by  [Se $07$, Theorem 5],
the bounds:

$$
v_\l(A_0)\leqslant
\left\{
\begin{array}{lll}

    2m + 2&  \iif & \l = 2\\
     2m&  \iif & \l > 2  \;\;\;{\rm and}\;\;  t = 1\;\;{\rm  or}\;\; 2\\
      m&  \iif &   t = 3, 4 \;\;\oor \;\; 6\\
      0 & \iif & t = 5  \;\;{\rm or}\;\;
t > 6.
\end{array}
\right.
$$

Since  $v_\l(A) =v_ \l(A_o)$  if  $\l > 2$  and  $v_\l(A) \leqslant v_\l(A_o) + 1$  if  $\l = 2$, we
obtain a bound for  $v_\l(A)$  which is  $\leqslant  M(k,\l)$.

\medskip

\ni {\it Remarks}. 1) Note the case  $\l = 2$, where the  $M(k,\l)$  bound  $2m+3$  can be
attained.

   2) In the case  $t = 6$, the bound  $v_\l(A_o) \leqslant  m$  given above can be
replaced by  $v_\l(A_o) = 0$, but this is not important for what we are doing here.

\subsection{The Del Pezzo case : degree $7$}

   This is a trivial case; 
there are $3$ exceptional curves on  $S$ (over  $\overline{k}$),
and only one of them meets the other two. It is thus stable under  $A$, and by
blowing it down, one is reduced to the degree $8$ case. [This case does
not occur if one insists, as in [DI 08], that the rank of  $\Pic(S)^A$  be equal to  $1$.]

\subsection{The Del Pezzo case : degree $6$}

   Here the surface  $S$  has $6$ exceptional curves (over  $\overline{k}$); their incidence
graph  $\Sigma$   is an hexagon. There is a natural homomorphism
     $$g:  \Aut(S) \rightarrow  \Aut(\Sigma)$$
and its kernel   $T$  is a $2$-dimensional torus. Put  $A_o = A  \cap  T(k)$.
The index of  $A_o$  in  $A$  is a divisor of $12$. By [Se $07$, Theorem4], we have

$$v_\l(A_o)\leqslant 
\left\{
\begin{array}{llll}

      2m  & \iif &  t = 1 \;\; \oor \;\; 2     &({\rm i.e.} \;\;\iif\;\;  \varphi(t) = 1)\\
        m  & \iif & t = 3, \ 4 \;\; \oor \;\;  6  &({\rm i.e}. \;\;\iif\;\;  \varphi(t) = 2)\\
      0    & \iif & t = 5  \;\;\oor \;\; t >6.&

\end{array}
\right.
$$

Hence:
$$v_\l(A)\leqslant
\left\{
\begin{array}{lll}

   2m+2&  \iif & \l = 2\\
   2m+1&  \iif & \l = 3\\
   2m&    \iif & \l > 3  \;\;{\rm and}\;\;  t = 1 \;\;\oor \;\; 2\\
   m&     \iif & t = 3,\  4 \;\;\oor \;\; 6\\
   0&     \iif & t = 5  \;\;\oor \;\;  t > 6.
\end{array}
\right.
$$
These bounds are $\leqslant M(k,\l)$.

\medskip

\ni {\it Remarks.} 1) Note the case  $t = 6$, where the bound  $m$  can actually
be attained.

  2) In the case  $t = 4$, the bound  $v_\l(A) \leqslant  m$  given above can be replaced
by  $v_\l(A) = 0$.

\subsection{The Del Pezzo case : degree $5$}

   As above, let  $\Sigma$  be the incidence graph of the exceptional curves of  $S$. Since
$\deg(S) \leqslant 5$, the natural map  $\Aut(S) \rightarrow \Aut(\Sigma)$  is injective. We can thus
identify  $A$  with its image  in  $\Aut(\Sigma)$. In the case  $\deg(S) = 5$, the graph $\Sigma$ is the Petersen graph, and $\Aut(\Sigma)$  is isomorphic to the symmetric group  $S_5$. This shows that

$$v_\l(A)\leqslant
\left\{
\begin{array}{lll}

   3&  \iif & \l = 2\\
     1&  \iif & \l = 3 \;\;\oor\;\; 5\\
     0&  \iif & \l > 5,
\end{array}
\right.
$$
and we conclude as before.

\subsection{The Del Pezzo case : degree $4$}

   This case is similar to the preceding one. Here  $\Aut(\Sigma)$  is isomorphic
to the group  $2^4.S_5 = \Weyl(D_5)$; its order is  $2^7.3.5$. We get the
same bounds as above, except for  $\l = 2$  where we find  $v_\l(A) \leqslant 7$,
which is  $\leqslant  M(k,2)$ [recall that  $M(k,2) = 2m+3$  and that  $m \geqslant 2$
for  $\l = 2$].

\subsection{The Del Pezzo case : degree $3$}

   Here   $S$  is a smooth cubic surface, and  $A$  embeds in $\Weyl(E_6)$, a group of order  $2^7.3^4.5$. This gives a
bound for  $v_\l(A)$  which gives what we want, except when $\l = 3$. 
In the case  $\l = 3$, it gives  $v_\l(A) \leqslant 4$, but Theorem $2.1$ claims  $v_\l(A) 
\leqslant 3$  unless  $k$  contains a primitive cubic root of unity. We 
thus have to prove the following lemma:

\medskip

\ni
{\bf Lemma $4.2$} - {\it Assume that  $|A| = 3^4$ , that  $A$  acts faithfully on a smooth
cubic surface  $S$  over  $k$, and that $\charb(k) \neq 3$. Then  $k$  contains a primitive cubic root of unity.}

\medskip

{\it Proof}. The structure of  $A$  is known since  $A$  is isomorphic to a $3$-Sylow
subgroup of $\Weyl(E_6)$. In particular the center  $Z(A)$   of  $A$  is cyclic
of order  $3$  and is contained in the commutator subgroup  of  $A$. Since  $A$
acts on $S$, it acts on the sections of the anticanonical sheaf of  $S$;
we get in this way a faithful linear representation  $r : A  \rightarrow  \GL_4(k)$.
Over  $\overline{k}$,  $r$  splits as  $r= r_1 + r_3$  where  $r_1$  is $1$-dimensional
and  $r_3$  is irreducible and $3$-dimensional. If  $z$  is a non trivial
element of  $Z(A)$, the eigenvalues of  $z$  are  $\{1,r,r,r\}$  where   $r$  is
a primitive third root of unity. This shows that  $r$  belongs to  $k$.

\subsection{The Del Pezzo case : degree $2$}

   Here  $A$  embeds in $\Weyl(E_7)$, a group of order  $2^{10}.3^4.5.7$. This gives
a bound for  $v_\l(A)$, but this bound is not good enough. However, the surface
$S$  is a $2$-sheeted covering of $\P_2$  (the map  $S  \rightarrow  \P_2$  being the
anticanonical map) and we get a homomorphism  $g: A \rightarrow \PGL_3(k)$
whose kernel has order $1$ or $2$. We then find the same bounds for  $v_\l(A)$  as in
\S 4.2, except that, for $\l = 2$, the bound is  $2m+2$  instead of  $2m+1$.

\subsection{The Del Pezzo case : degree 1}

   We use the linear series  $|-2 K_S|$.
It gives a map  $g: S \rightarrow \P_3$
whose image is a quadratic cone  $Q$, cf. e.g. [De $80$, p.68]. This
realizes  $S$  as a quadratic covering of  $Q$. If  $B$  denotes the
automorphism group of  $Q$  defined by  $A$, we have  $v_\l(A) = v_\l(B)$  if $\l > 2$
and  $v_\l(A) \leqslant v_\l(B) + 1$  if  $\l = 2$. But  $B$  is isomorphic to a subgroup
of  $k^* \times \Aut(C)$, where  $C$  is a curve of genus  $0$. This implies

$$v_\l(B)\leqslant 
\left\{\begin{array}{lll}

                     m + m + 1&    {\rm if}&  \l = 2\\
                m + m&        {\rm if}&  t = 1\\
                       0 + m&        {\rm if}&  t = 2,\; \l > 2\\
                       0 + 0&        {\rm if} & t > 2.
\end{array}
\right.$$

 The corresponding bound for  $v_\l(A)$  is $\leqslant  M(k,\l)$.

\medskip

\ni This concludes the proof of Theorem 4.1 and hence of Theorem 2.1.
   
   \bigskip
 \section{Structure and  conjugacy properties of 
       $\l$-subgroups of $\Cr(k)$}
   
   \medskip
\subsection{The $\l$-subgroups of $\Cr(k)$}

\medskip
The main theorem (Theorem 2.1) only gives information on the order of an $\l$-subgroup $A$ of $\Cr(k)$, assuming as usual that $\l \neq \charb(k)$. As for the structure of $A$, we have:

\medskip
\ni
{\bf Theorem $5.1$}. (i) {\it If $\l > 3$, $A$ is abelian of rank $\leqslant 2$ }(i.e. can be
generated by two elements).

 {\rm (ii)}. {\it If $\l = 3$ }(resp. $\l = 2$) {\it $A$ contains an abelian normal subgroup of
 rank $\leqslant 2$ with index $\leqslant 3$ } (resp. with index $\leqslant 8$).

 \medskip
 \ni {\it Proof}. Most of this is a consequence of the results of [DI 07]; see also [Bl 06] and [Be07]. The only case which does not seem to be explicitly in [DI 07] is the case $\l = 2$, when
 $A$ is contained in $\Aut(S)$, where $S$ is a conic bundle. Suppose we are in that case and let $f: S \rightarrow C$  and  $A_o$, $B$  be as in ¤4.3, so that we have an exact sequence            $ 1 \rightarrow B \rightarrow A \rightarrow A_o \rightarrow 1$,                  with $A_o \subset \Aut(C)$, and $B \subset \Aut(F)$ where  $F$ is the generic fiber of $f$ (which is a genus zero curve over the function field $k(C)$ of $C$). We use the following lemma:
 
 \medskip
 
 \ni {\bf Lemma $5.2$.} {\it Let $a\in A$ and $b\in B$ be such that $a$ normalizes the cyclic group $\langle b \rangle$  generated by $b$. Then $aba^{-1}$ is equal to $b$ or to $b^{-1}$.}
 
 \medskip
 \ni {\it Proof of the lemma}. Let $n$ be the order of $b$. If $n = 1$  or $2$, there is nothing to prove. Assume $n > 2$. By extending scalars, we may also assume that $k$ contains the primitive $n$-th roots of unity. Since $b$ is an automorphism of $F$ of order $n$, it fixes two rational points of $F$ which one can distinguish by the eigenvalue of $b$ on their tangent space: one of them gives a primitive $n$-th root of unity $z$, and the other one gives $z'=z^{-1}$. [Equivalently, $b$ fixes two sections of $f:  S \rightarrow C$.] The pair $(z,z')$ is canonically associated with $b$. Hence the pair associated with $aba^{-1}$ is also $(z,z')$. On the other hand, if $aba^{-1} = b^i$ with $i\in\textbf{Z}/n\textbf{Z}$, then the pair associated to $a^i$ is $(z^i,z'^i)$. This shows that $z^i$ is equal to either $z$ or $z^{-1}$, hence $i \equiv 1$ or $-1$  (mod  $n$). The result follows.
 \medskip
 
 \ni {\it End of the proof of Theorem 5.1 in the case $\l = 2$.} Since $B$ is a finite $2$-subgroup of a $k(C)$-form of $\PGL_2$, it is either cyclic or dihedral. In both cases, it contains a characteristic subgroup $B_1$ of index $1$ or $2$  which is cyclic. Similarly, $A$ has a cyclic subgroup $A_1$ which is of index $1$ or $2$. Let $a\in A$ be such that its image in  $A_o$ generates $A_1$. If $b$ is a generator of  $B_1$, Lemma $5.2$ shows that $a^2$
 commutes with $b$. Let $\langle b,a^2 \rangle$ be the abelian subgroup of $A$ generated by $b$ and $a^2$. It is normal in  $A$, and the inclusions  $\langle b,a^2 \rangle  \subset     \langle b,a \rangle   \subset   B.\langle a\rangle \subset  A$    show that its index in $A$ is at most $8$.
 
\medskip
\ni {\it Remark}. Similar arguments can be applied to prove a Jordan-style result on the finite subgroups of $\Cr(k)$, namely:

\medskip
\ni{\bf Theorem $5.3$.} {\it There exists an integer $J > 1$, independent of the field $k$, such that every finite subgroup $G$ of $\Cr(k)$, of order prime to $\charb(k)$, contains an abelian normal subgroup $A$ of rank $\leqslant 2$, whose index in $G$ divides $J$.}

\medskip
The proof follows the same pattern: the conic bundle case is handled via Lemma $5.2$ and the Del Pezzo case via the fact that $G$ has a subgroup of bounded index which is contained in a reductive group of rank $\leqslant 2$, so that one can apply the usual form of Jordan's theorem to that group. As for the value of  $J$, a crude computation shows that one can take  $J = 2^{10}.3^4.5^2.7$; the exponents of $2$ and $3$ can be somewhat lowered, but those of $5$ and $7$ cannot since $\Cr(\C)$  contains $A_5 \times A_5$ 
and  $\PSL_2(\F_7)$.
\bigskip
\subsection{The cases t = $3,4,6$}
\medskip
  More precise results on the structure of $A$ depend on the value of the invariant $t = t(k,\l)$. Recall that $t = 1,2,3,4$ or $6$ if $A \neq 1$, cf. Cor.2.2. We shall only consider the cases $t = 3,4$  or $6$ which are the easiest. See [DI 08, \S4] for a (more difficult) conjugation theorem which applies when $t = 1$ or $2$.
   Recall (cf. $\S3.2$)  that  $A$  is  said to be {\it toral} if there exists a $2$-dimensional subtorus $T$ of $\Cr$ (in the sense of [De 70]) such that $A$ is contained in $T(k)$. We have:
   
   \medskip
   \ni
   {\bf Theorem $5.4$}.\ {\it  Assume that $t = 3,4$ or $6$. Then}:
   
  \ni (a) $A$ {\it is cyclic of order $\l^n$ with} $n \leqslant m$.
   
  \ni (b) $A$ {\it is toral, except possibly if} $|A| = 5$.
  
  \ni (c) {\it If $A'$ is a subgroup of $\Cr(k)$ of the same order as $A$, then $A'$ is conjugate to $A$ in $\Cr(k)$, except possibly if} $|A| = 5$.
  \medskip
  
  Note that the hypothesis $t = 3,4$ or $6$ implies $\l \geqslant 5$. Moreover, if $\l = 5$, then $t=4$ and, if $\l=7$, then $t=3$ or $6$.
  \medskip
  
  \ni {\it Proof of} (a) {\it and} (b). We follow the same method as above, i.e. we view $A$ as a subgroup of $\Aut(S)$, where $S$ is either a conic bundle or a Del Pezzo surface. The bounds given in \S4.3 show that $A = 1$ if $S$ is a conic bundle (this is why this case is easier than the case $t = 1$ or $2$). Hence we may assume that $S$ is a Del Pezzo surface. Let $d$ be its degree. We have an exact sequence:
  
    $1\ \rightarrow \ G(k) \ \rightarrow \ \ \Aut(S) \ \rightarrow \ E \ \rightarrow \ 1$,
    
\ni   where $G = \Aut(S)^o$ is a connected linear group of rank $\leqslant 2$ and $E$ is a subgroup of a Weyl group $W$ depending on $d$ (e.g. $W =$ Weyl$(E_8)$ if $d=1$).

  Consider first the case $\l > 7$. The order of $W$ is not divisible by $\l$ ; hence $A$ is contained in $G(k)$. Since $A$ is commutative, there exists a maximal torus  $T$ of $G$ such that $A$ is contained in the normalizer $N$ of $T$, cf. e.g. [Se 07, \S3.3]; since
 $\l > 3$, the order of $N/T$ is prime to $\l$, hence $A$ is contained in $T(k)$ and this
 implies dim$(T) \geqslant 2$ by [Se 07, \S4.1]. This proves (b),  and (a) follows from Lemma 5.5 below.
 
   Suppose now that $\l = 5$ or $7$, and let $n = v_\l(A)$. If $n = 1$ and $\l = 5$ , there is nothing to prove. If $n = 1$ and $\l = 7$, then (a) is obvious and (b) is proved in [DI 08, prop.3] (indeed Dolgachev and Iskovskikh prove (b) when $v_\l(A) = 1$, and they also prove (c) for $\l = 7$). We may thus assume that $n > 1$. If $d \leqslant 5$, then $G = 1$ and 
$A$ embeds in $E$; but $E$ does not contain any subgroup of order $\l^2$ (see the tables in [DI 07] and [Bl 06]); hence this case does not occur. If $d > 5$, then the order of $E$ is prime to $\l$, hence $A$ is contained in $G(k)$ and the proof above applies.

\medskip
 \ni {\it Proof of} (c). By (b), we have  $A \ \subset  \ T(k)$ and $A'\ \subset \ T'(k)$ where $T$ and $T'$ are 2-dimensional subtori of $\Cr$. By Lemma 5.5 below, these tori are isomorphic; by a standard argument (see e.g. [De 70, \S6] this implies that $T$ and $T'$ are conjugate by an element of $\Cr(k)$; moreover $A$ (resp. $A'$) is the unique subgroup of order $\l^n$ of $T(k)$ (resp. of $T'(k)$). Hence $A$ and $A'$ are conjugate in $\Cr(k)$.
 
 \medskip
 \ni {\it Remark.} The case $|A| = 5$ is indeed exceptional: there are examples of such 
 $A$'s which are not toral, cf. [Be 07], [Bl 06], [DI 07].
 \medskip
 \subsection {A uniqueness result for $2$-dimensional tori}
 
 We keep the assumption that $t = 3,4$ or $6$. We have seen in \S3.2.2 that there exists a 2-dimensional $k$-torus  $T$  such that $T(k)$ contains an element of order $ \l$.
 \medskip
 
 \ni {\bf Lemma} $5.5.$ (a) {\it Such a torus is unique, up to $k$-isomorphism.}
 
 (b) {\it If $n \leqslant m = m(k,\l)$, then $T(k)[\l^n]$ is cyclic of order $\l^n$.}
 
\medskip
\ni {\it Proof of} (a). Let  $L =  \Hom_{k_s}(\G_m,T)$ be the group of cocharacters of $T$. It is a free \Z-module of rank $2$, with an action of $\Gamma_k = \Gal(k_s/k)$. If we identify $L$ with $\Z^2$, this action gives a homomorphism $r: \Gamma_k \rightarrow \GL_2(\Z)$ which is well defined up to conjugation. Let  $G$  be the image of $r$. Since $G$ is a finite subgroup of $\GL_2(\Z)$, its order divides $24$, and hence is prime to $\l.$

  The $\Gamma_k$-module $T(k_s)[\l]$ of the $\l$-division points of $T(k_s)$ is canonically isomorphic to $L/\l L \otimes \mu_\l$, where $\mu_\l$ is the group of $\l$-th roots of unity in $k_s$. This shows that $L/\l L$ contains a rank-1 submodule $I$ which is isomorphic to the dual $\mu_\l^*$ of $\mu_\l$. The action of $G$ on $L/\l L$ is semisimple since $|G|$ is prime to $\l$. Hence there exists a rank-1 submodule $J$ of $L/\l L$ such that $L/\l L = I \bigoplus J$. By a well-known lemma of Minkowski (see e.g. [Se 07, Lemma 1]), the action of  $G$  on  $L/\l L$ is faithful. This shows that $G$ is commutative. Moreover, the character giving the action of $\Gamma_k$ on $I$ has an image which is cyclic of order $t$. Since $t = 3,4$ or $6$, this shows that $G$ contains an element of order $3$ or $4$. One checks that these properties imply $G \subset  \SL_2(\Z)$  i.e. $\det(r) = 1$, hence the $\Gamma_k$-modules $I$ and $J$ are dual of each other, i.e.  $J \simeq \mu_\l$. We thus have $L/\l L \simeq \mu_l \oplus \mu_\l^*$ . We may then identify $r$ with the homomorphism $\Gamma_k\ \rightarrow\ C_t\ \rightarrow \ \GL_2(\Z)$, where  $C_t$ is the Galois group of $k(\mu_\l/k)$ and  $C_t \ \rightarrow \ \GL_2(\Z)$ is an inclusion. Since any two such inclusions only differ by an inner automorphism of $\GL_2(\Z)$, this shows that the $\Gamma_k$-module $L$ is unique, up to isomorphism; hence the same is true for $T$.
  
  \ni {\it Proof of} (b).  Assertion (b) follows from the description of $T$ given in $\S3.2.2$. It can also be checked by writing explicitly the $\Gamma_k$-module $L/\l^n L$; when $n \leqslant m$ this module is isomorphic to the direct sum of $\mu_{\l^n}$ and its dual.
  
 \medskip
 \ni {\it Remarks.}
 
 1). If $n > m$ we have  $T(k)[\l^n] = T(k)[\l^m].$ This can be seen, either by a direct computation of $\l$-adic representations, or by looking at $\S3.2.2.$ 
 
 2) When $t = 1$ or $2$, it is natural to ask for a 2-dimensional torus $T$ such that $T(k)$ contains  $\Z/\l Z \oplus \Z/\l \Z$. Such a torus exists, as we have seen in $\S3.2$. If $\l > 2$,  it is unique, up to isomorphism. There is a similar result for $\l = 2$, if one asks not merely that $T(k)$ contains $\Z/2\Z \times \Z/2\Z$ but that it contains $\Z/4\Z \times \Z/4\Z$.
 
 \bigskip
 
 \section{The Cremona groups of rank $>  2$}
 For any $r > 0$ the Cremona group $\Cr_r(k)$ of rank  $r$ is defined  as the
 group $\Aut_k k(T_1,...,T_r)$  where $(T_1,...,T_r)$  are $ r$ indeterminates.  When $r > 2$ not much seems to be known on the finite subgroups of  $\Cr_r(k)$, even in the classical case  $k = \C$. For instance:
 
 \smallskip
 
 \ni $6.0$.  {\it Does there exist a finite group which is not embeddable in $\Cr_3(\C)$} ?
 
 This looks very likely. It is natural to ask for much more, for instance :

 \medskip
 \ni $6.1$ (Jordan bound, cf. Theorem 5.3). {\it Does there exist an integer $N(r) > 0$ , depending only on $r$, such that, for every finite subgroup  $G$  of  $\Cr_r(k)$
 of  order prime to $\charb(k)$, there exists an abelian normal subgroup $A$ of  $G$, of rank  $ \leqslant r$,  whose index divides $N(r)$ }?
 
 Note that this would imply that, for $ \l$ large enough (depending on $r$),
 every finite $\l$-subgroup of $\Cr_r(k)$ is abelian of rank $ \leqslant r$.
 
 \smallskip
 
 \ni $6.2$ (cf. [Se 07, $\S6.9$] ). {\it Is it true that  $r \geqslant \varphi(t)$ if  $\Cr_r(k)$ contains an element of order $\l$ }?
 
 \smallskip

 \ni $6.3$. {\it Let $ G \subset \Cr_r(k)$ be as in $6.1$, and assume that $k$ is small} (cf.  $\S2.3$).{\it  Is it true that $|G|$ is bounded by a constant depending only on $r$ and the cyclotomic invariants $(t,m)$ of $k$} ?
 
 \smallskip
 
   If the answer to $6.3$ is `` yes '' we may define $M_r(k)$  as the l.c.m. of all such $ |G|$'s, and  ask for an estimate of  $M_r(k)$. For instance, in the case  $r = 3$ :
   
   \smallskip
   
\smallskip   \ni $6.4$. {\it Is it true that  $M_3(k)$ is equal to  $M_1(k)M_2(k)$} ?

\smallskip

   If $k$ is finite with $q$ elements, this means (cf.  $\S2.5$):
   
   \smallskip
   
   \ni $6.5.$ {\it Is it true that
   
   \medskip
   
   $M_3(k) = \left\{  
   \begin{array}{ll}
        3.(q^2 - 1)(q^4 - 1)(q^6 - 1)& {\it if} \  q \equiv 4 \ {\it or} \; 7 \pmod 9\\
          (q^2 - 1)(q^4 - 1)(q^6 - 1) &  {\it otherwise}  \ ?
     \end{array}
     \right.
     $
     }
     
     \medskip
For larger  $r$'s the polynomial $(X^2  - 1)(X^4 - 1)(X^6 - 1)$ of $6.5$ should be replaced by the polynomial $P_r(X)$ defined by the formula

\smallskip

    \ \ $P_r(X) = \displaystyle{\prod_{d} \Phi_d(X)^{[r/\varphi(d)]}}$  \ ,

\ni where  $\Phi_d(X)$ is the $d$-th cyclotomic polynomial.
\medskip

\ni {\it Examples}. $P_4(X) = (X^6 - 1)(X^8 - 1)(X^{10} - 1)(X^{12} - 1);
P_5(X) = (X^2 - 1)P_4(X).$

\medskip

 With this notation, the natural question to ask  seems to be :

\smallskip

\ni $6.6$. {\it Is it true that there exists an integer  $c(r) > 0$  such that $M_r(\F_q)$ divides
$c(r).P_r(q)$  for every  $q$ ?}
\smallskip

  Unfortunately, I do not see how to attack these questions; the method used for rank 2 is based on the detailed knowledge of the ``minimal models'', and this is not available for higher ranks.

  \bigskip
  
 \ni {\it Acknowledgment.} I wish to thank A.Beauville for a series of e-mails in $2003-2005$ which helped me to correct the naive ideas I had on the Cremona groups.

\begin{center}
       {\bf References}
\end{center}

\bigskip

[Be 07] A.Beauville, {\it p-elementary subgroups of the Cremona group},  J. Algebra
     $\bf{314}$ (2007), 553-564.

[Bl 06]  J.Blanc, {\it Finite abelian subgroups of the Cremona group of the
     plane}, Univ.Gen\`eve, th\`ese \no $3777$ (2006). See also C.R.A.S. $\bf{344}$
    ($2006$), 21-26.

[De 70]  M.Demazure, {\it Sous-groupes alg\'ebriques de rang maximum du groupe de
     Cremona}, Ann.Sci.ENS ($4$) $\bf{3}$ ($1970$), 507-588.

[De 80]  M.Demazure, {\it Surfaces de Del Pezzo}, I-IV, Lect.Notes in Math. $\bf{777},$
    Springer-Verlag, $1980$, 21-69.

[Do 07]  I.V.Dolgachev, {\it Topics in Classical Algebraic Geometry, Part I},
     Lecture notes, Univ. Michigan, Ann Arbor, $2007$.

[DI 07]  I.V.Dolgachev and V.A.Iskovskikh, {\it Finite subgroups of the plane
    Cremona group}, ArXiv:math/0610595v2, to appear in {\it Algebra, Arithmetic
    and Geometry, Manin's Festschrift}, Progress in Math. Birkh\"auser Boston,
    $2008$.

[DI 08]  I.V.Dolgachev and V.A.Iskovskikh, {\it On elements of prime order in
    the plane Cremona group over a perfect field}, ArXiv:math/0707.4305,
    to appear.

[Is 79]  V.A.Iskovskikh, {\it Minimal models of rational surfaces over arbitrary
    fields} (in Russian), Izv.Akad.Nauk $\bf{43}$ ($1979$) , 19-43; English translation:
    Math. USSR Izvestija $\bf{14}$ (1980), 17-39. 

[Is 96]  V.A.Iskovskikh, {\it Factorization of birational maps of rational
    surfaces from the viewpoint of Mori theory} (in Russian), Uspekhi
    Math.Nauk $\bf{51}$(1996), 3-72; English translation: Russian Math.Surveys
    $\bf{51} $ (1996), 585-652. 

[Ko 96]  J.Koll\'ar, {\it Rational Curves on Algebraic Varieties}, Ergebn.Math. ($3$)
    $\bf{32}$, Springer-Verlag, 1996. 

[Ma 66]  Y.I.Manin, {\it Rational surfaces over perfect fields} (in Russian, with
    English r\'esum\'e), Publ.Math.IHES $\bf{30} $ (1966), 415-475.

[Ma 86]  Y.I.Manin, {\it Cubic Forms}: {\it Algebra, Geometry, Arithmetic}, 2nd edition,
     North Holland, Amsterdam, 1986.

[Se 07]  J-P.Serre, {\it Bounds for the orders of the finite subgroups of  $G(k)$,}
    in {\it Group Representation Theory}, eds. M.Geck, D.Testerman \& J. Th\'evenaz,
    EPFL Press, Lausanne, 2007, 403-450.
    
[Se 08]  J-P.Serre, {\it Le groupe de Cremona et ses sous-groupes finis}, S\'em.
Bourbaki 2008/2009, expos\'e 1000.

[Vo 98]  V.E.Voskresenski\u{\i}, {\it Algebraic Groups and Their Birational
     Invariants}, Translations Math. Monographs $\bf{179}$, AMS, 1998. 

\bigskip

\ni
Coll\`ege de France

\ni
3, rue d'Ulm

\ni
F-$75231$ Paris Cedex $05$
\medskip

\ni
e-mail: serre@noos.fr

\end{document}